\documentclass[final]{elsarticle}

\usepackage{amssymb}
\usepackage{amsfonts}
\usepackage{amsmath}
\usepackage[english,francais]{babel}
\usepackage[latin1]{inputenc}

\newtheorem{theorem}{Theorem}[section]

\newtheorem{e-proposition}[theorem]{Proposition}

\newtheorem{e-definition}[theorem]{Definition\rm}

\def\og{\leavevmode\raise.3ex\hbox{$\scriptscriptstyle\langle\!\langle$~}}
\def\fg{\leavevmode\raise.3ex\hbox{~$\!\scriptscriptstyle\,\rangle\!\rangle$}}

\newcommand{\RR}{\field{R}}

\newcommand{\be}{\begin{equation}}
\newcommand{\ee}{\end{equation}}
\newcommand{\ba}{\begin{array}}
\newcommand{\ea}{\end{array}}
\newcommand{\bea}{\begin{eqnarray}}
\newcommand{\eea}{\end{eqnarray}}
\newcommand{\bee}{\begin{eqnarray*}}
\newcommand{\eee}{\end{eqnarray*}}

\catcode`@=11
\renewcommand\appendix{\bigskip {\noindent\Large \bf Appendix}\par
  \setcounter{section}{0}%
  \setcounter{subsection}{0}%
  \renewcommand\thesection{\@Alph\c@section}}
\catcode`@=12

\newenvironment{acknowledgement}{\noindent{\bf Acknowledgement.~}}{}

\def\RR{\mathbb{R}}
\def\R{{\mathbb R}}

\def\lim{\mathop{\rm lim}}

\def\log{{\rm log}}

\def\lsl{\frac{\lambda_s}{\lambda}}

\def\fref#1{{\rm (\ref{#1})}}

\def\calE{{\mathcal E}}

\def\pa{\partial}

\begin{document}

\begin{frontmatter}
\selectlanguage{english}
{
\title{Blow up dynamics for smooth data equivariant solutions to the energy critical Schr\"odinger map problem}
}



\selectlanguage{english}
\author[authorlabel1]{Frank Merle}
\ead{merle@math.u-cergy.fr}
\author[authorlabel2]{Pierre Rapha\"el}
\ead{pierre.raphael@math.univ-toulouse.fr}
\author[authorlabel3]{Igor Rodnianski}
\ead{irod@math.princeton.edu
}

\address[authorlabel1]{Universit{\'e} de Cergy Pontoise et IHES,
2 av Adolphe Chauvin, 95 302 Cergy Pontoise, France}

\address[authorlabel2]{  Institut de Math\'ematiques, Universit\'e Paul Sabatier,
118 Route de Narbonne, 31062 Toulouse Cedex, France}

\address[authorlabel3]{Department of Mathematics, Princeton University,
Fine Hall, Washington Road, NJ 08544-1000, USA}

\begin{abstract}
We consider the energy critical Schr\"odinger map $\pa_tu=u\wedge \Delta u$ to the 2-sphere for equivariant initial data of homotopy number $k=1$. We show the existence of a set of smooth initial data arbitrarily close to the ground state harmonic map $Q_1$ in the scale invariant norm $\dot{H}^1$ which generates finite time blow up solutions. We give in addition a sharp description of the corresponding singularity formation which occurs by concentration of a universal bubble of energy  $$u(t,x)-e^{\Theta^* R}Q_1\left(\frac{x}{\lambda(t)}\right)\to u^*\ \ \mbox{in} \ \ \dot{H}^1 \ \ \mbox{as} \ \ t\to T$$ 
where $\Theta^* \in \R$, $ u^* \in  \dot{H}^1$, $R$ is a rotation and the concentration rate is given for some $\kappa(u)>0$ by $$\lambda(t)=\kappa(u)\frac{T-t}{|\log (T-t)|^2}(1+o(1)) \ \ \mbox{as} \ \ t\to T.$$ Full details of the proofs will appear in the companion
paper \cite{MRR}.
\vskip 0.5\baselineskip

\selectlanguage{francais}
\noindent{\bf R\'esum\'e}
\vskip 0.5\baselineskip
\noindent

Nous consid\'erons l'application de Schr\"odinger sur la 2-sph\`ere \'energie critique $\pa_tu=u\wedge \Delta u$ pour des donn\'ees initiales \`a sym\'etrie \'equivariante et de degr\'e $k=1$. Nous exhibons un ensemble de donn\'ees initiales r\'eguli\`eres arbitrairement proches dans la topologie invariante d'\'echelle $\dot{H}^1$ de l'application harmonique d'\'energie minimale $Q_1$ qui engendrent des solutions explosives en temps fini. Nous donnons une description fine de la formation de singularit\'e qui correspond \`a la concentration d'une bulle universelle  d'\'energie $$u(t,x)-e^{\Theta^*R}Q_1\left(\frac{x}{\lambda(t)}\right)\to u^*\ \ \mbox{in} \ \ \dot{H}^1$$ 
o\`u $\Theta^* \in \R$, $ u^* \in  \dot{H}^1$,  $R$ est une rotation et la vitesse de concentration est donn\'ee pour une certain $\kappa(u)>0$ par :  $$\lambda(t)=\kappa(u)\frac{T-t}{|\log (T-t)|^2}(1+o(1)) \ \ \mbox{quand} \ \ t\to T.$$
\end{abstract}

\end{frontmatter}

\selectlanguage{francais}
\section*{Version fran\c{c}aise abr\'eg\'ee}

Nous consid\'erons l'application de Schr\"odinger  \'energie critique sur la 2-sph\`ere 
\be
\label{nlsmapfracnsi}
\left\{\begin{array}{ll}\pa_tu=u\wedge \Delta u,\\u_{|t=0}=u_0 \in \dot{H}^1 \end{array}\right . \ \ (t,x)\in \Bbb R\times \Bbb R^2, \  \ u(t,x)\in \Bbb S^2.
\ee Ce syst\`eme appartient \`a une classe d'\'equations g\'eom\'etriques \'energie critique qui inclut le flot  parabolique de la chaleur harmonique cf. \cite{S}, \cite{Q}, \cite{QT},  \cite{PT}, \cite{heatflow} et les applications de type onde, cf \cite{RaphRod},  et appara\^it notamment en ferromagn\'etisme en relation avec les \'equations de Landau-Lifschitz. Ce syst\`eme est Hamiltonien et le flot laisse invariante l'\'energie de Dirichlet  
\be
\label{dircehitkfrac}
E(u(t))=\int_{\Bbb R^2}|\nabla u(t,x)|^2dx=E(u_0).
\ee
Nous consid\'erons des flots \`a sym\'etrie k-\'equivariante $$u(t,x)=e^{k\theta R}\left|\begin{array}{lll} u_1(t,r)\\ u_2(t,r)\\u_3(t,r),\end{array}\right .\ \ R=\left(\begin{array}{lll} 0 & -1 & 0 \\ 1 & 0 &0\\ 0 & 0& 0
\end{array}
\right ),$$
 o\`u $(r,\theta)$ d\'esignent les coordonn\'ees polaires sur $\Bbb R^2$, et o\`u $k\in \Bbb Z^*$ est le degr\'e de l'application. Dans ce cas, le probl\`eme de Cauchy est bien pos\'e localement en temps dans $\dot{H}^1$ par \cite{CSU}, \cite{T3}, \cite{T4}. Pour une donn\'ee initiale g\'en\'erale le resultat est connu uniquement pour donn\'ee petite \cite{BT2}.

 Le minimiseur de l'\'energie de Dirichlet \fref{dircehitkfrac} \`a degr\'e fix\'e est explicitement donn\'e par $$Q_k(x)=e^{k\theta R}\left|\begin{array}{lll}\frac{2r^k}{1+r^{2k}} \\ 0\\\frac{1-r^{2k}}{1+r^{2k}}\end{array}\right .$$ et engendre une solution stationnaire de \fref{nlsmapfracnsi}. Pour $k\geq 3$, cette solution est stable et m\^eme asymptotiquement stable d'apr\`es Gustafson, Nakanishi, Tsai \cite{T2}. Pour $k=1$, le premier r\'esultat d'instabilit\'e de $Q\equiv Q_1$ dans la topologie invariante d'\'echelle $\dot{H}^1$ est donn\'e par Bejenaru et Tataru \cite{BT}. Dans la continuation des travaux sur l'\'equation de Schr\"odinger $L^2$ critique \cite{MR1}, \cite{MR2}, \cite{MR3}, \cite{MR4}, \cite{MR5}, \cite{R1} et l'application d'ondes \'energie critique sur la 2-sph\`ere \cite{RaphRod}, nous obtenons le premier r\'esultat d'explosion en temps fini:\\
 
 {\bf Th\'eor\`eme} [Explosion pour $k=1$] {\it Il existe un ensemble de donn\'ee initiales r\'eguli\`eres \`a sym\'etrie \'equivariante, de degr\'e $k=1$ et arbitrairement proches de $Q_1$ dans $\dot{H}^1$ telles que la solution correspondante de \fref{nlsmapfracnsi} explose en temps fini $T<+\infty$ par concentration d'une bulle universelle d'\'energie $$u(t,x)-e^{\Theta^* R}Q_1\left(\frac{x}{\lambda(t)}\right)\to u^*\ \ \mbox{dans} \ \ \dot{H}^1 \ \ \mbox{quand} \ \ t\to T$$ o\`u $\Theta^* \in \R$, $ u^* \in  \dot{H}^1$, et la vitesse de concentration est donn\'ee pour un certain $\kappa(u)>0$ par:  $$\lambda(t)=\kappa(u)\frac{T-t}{|\log (T-t)|^2}(1+o(1)) \ \ \mbox{quand} \ \ t\to T.$$}

Cette note est une version abr\'eg\'ee de \cite{MRR}.

\selectlanguage{english}


\section{Setting of the problem and main result}


In this paper we consider the energy critical Schr\"odinger map 
\be
\label{nlsmap}
\left\{\begin{array}{ll}\pa_tu=u\wedge \Delta u,\\u_{|t=0}=u_0 \in \dot{H}^1 \end{array}\right . \ \ (t,x)\in \Bbb R\times \Bbb R^2, \  \ u(t,x)\in \Bbb S^2.
\ee
This equation is related to the Landau-Lifschitz equation for ferromagnetism, and it 
belongs to a class of geometric evolution equations \cite{S}, \cite{Q}, \cite{QT},  \cite{PT}, \cite{heatflow}, including wave  maps and the harmonic heat flow, which have attracted a considerable attention in the past ten years. 
The Hamiltonian structure of the problem implies conservation of the Dirichlet energy \be
\label{dircehit}
E(u(t))=\int_{\Bbb R^2}|\nabla u(t,x)|^2dx=E(u_0)
\ee
which moreover is left unchanged by the scaling symmetry of the problem $u(t,x)\mapsto u_{\lambda}(t,x)=u(\lambda^2 t,\lambda x)$. 
The question of the global existence of all solutions or on the contrary the possibility of a finite blow up and singularity formation corresponding to a concentration of energy has been addressed recently in detail for the wave map problem -- the wave analogue of \fref{nlsmap} -- and the Yang-Mills equations, see \cite{T}, \cite{ST}, \cite{KS} for the large data wave map global regularity problem; 
\cite{RaphRod} and references therein, \cite{RSt}, \cite{KST} (see also \cite{S}, \cite{Q} \cite{QT}  \cite{PT}, \cite{GS}, \cite{heatflow} for the heat flow), and has been until now open for the Schr\"odinger map problem.\\

We shall focus on the case of solutions with k-equivariant symmetry
$$u(t,x)=e^{k\theta R}\left|\begin{array}{lll} u_1(t,r)\\ u_2(t,r)\\u_3(t,r),\end{array}\right .\ \ R=\left(\begin{array}{lll} 0 & -1 & 0 \\ 1 & 0 &0\\ 0 & 0& 0
\end{array}
\right ),$$
where $(r,\theta)$ are the polar coordinates on $\Bbb R^2$, and $k\in \Bbb Z^*$ is the homotopy number. 
In this case, the Cauchy problem is well-posed in $\dot{H}^1$ if the energy $E$ is sufficiently small, \cite{CSU} or, more 
generally, if the energy $E$ is sufficiently close to the minimum in a given homotopy class $k$, realized on a harmonic 
map $Q_k: \Bbb R^2\to \Bbb S^2$, \cite{T3} \cite{T4}. In the general case without symmetry 
the small energy data result is shown in (\cite{BT2}) and a conditional result for solutions with energy below that of the ground state
$Q_1$ is given in \cite{PS}.
In a given homotopy class, the minimizer of the Dirichlet energy \fref{dircehit} is explicitely given by the harmonic map $$Q_k(x)=e^{k\theta R}\left|\begin{array}{lll}\frac{2r^k}{1+r^{2k}} \\ 0\\\frac{1-r^{2k}}{1+r^{2k}}\end{array}\right .$$ which generates a stationary solution to \fref{nlsmap}. For large degree $k\geq 3$, this solution is stable, in fact, asymptotically stable by the result of Gustaffson, Nakanishi and Tsai \cite{T2}. For $k=1$ which corresponds to least energy maps, Bejenaru and Tataru \cite{BT} exhibit some instability mechanism of $Q\equiv Q_1$ in the scale invariant space $\dot{H}^1$.\\

 In the companion paper \cite{MRR} we give full details of the proof of the following result on formation of singularities 
 for the Schr\"odinger maps 
 arising from a set of smooth data arbitrarily close to $Q$ for $k=1$.  This 
 continues a series of works on the $L^2$ critical nonlinear Schr\"odinger equation \cite{MR1}, \cite{MR2}, \cite{MR3}, \cite{MR4}, \cite{MR5}, \cite{R1},  and the sharp description of a stable blow up  for the wave map problem to the 2-sphere \cite{RaphRod}.\\

\begin{theorem} [Existence and sharp description of a blow up regime for $k=1$]
\label{thmmain}
Let $k=1$. There exists a set of smooth  equivariant initial data of degree $k=1$ arbitrarily close to the ground state $Q_1$ in the $\dot{H}^1$ topology such that the corresponding solution to \fref{nlsmap} blows up in finite time through the concentration of a universal bubble of energy 
$$u(t,x)-e^{\Theta^* R}Q_1\left(\frac{x}{\lambda(t)}\right)\to u^*\ \ \mbox{in} \ \ \dot{H}^1\ \ \mbox{as} \ \ t\to T.$$ 
for some $\Theta^* \in \Bbb R$,  and at a speed given for some $\kappa(u)>0$ by: 
\be
\label{cneoheoe}
\lambda(t)=\kappa(u)\frac{T-t}{|\log (T-t)|^2}(1+o(1)) \ \ \mbox{as} \ \ t\to T
\ee
\end{theorem}

The blow up rate \fref{cneoheoe} is a natural candidate for a stable singularity formation, 
see \cite{heatflow} for the corresponding parabolic problem. However, in the Schr\"odinger map problem, as a consequence of 
a new instability mechanism,  \fref{cneoheoe} appears to represent a codimension one phenomena.


\section{Strategy of the proof}


{\bf step 1} Choice of gauge.\\
We describe the flow \fref{nlsmap} in the renormalized Frenet basis associated to the harmonic map $Q_1$: $$(e_r,e_\tau,Q_1), \ \ e_r=\frac{\pa_rQ_1}{|\pa_rQ_1|}, \ \ e_\tau=\frac{\pa_\tau Q_1}{|\pa_{\tau}Q_1|}, \ \ \pa_{\tau}=\frac 1r\pa_\theta.$$ We renormalise the map $$u(t,x)=e^{\Theta(t)R}v(s,y), \ \ \frac{ds}{dt}=\frac{1}{\lambda^2}, \ \ y=\frac{x}{\lambda}$$ and rewrite the equation for $w$ in the Frenet basis: $$v(s,y)=\alpha(s,y)e_r+\beta(s,y)e_\tau+(1+\gamma(s,y))Q_1,  \ \ \alpha^2+\beta^2+(1+\gamma)^2=1.$$ To leading order, the flow near $Q$ becomes a quasilinear Schr\"odinger equation:
\be
\label{renormalizeequation}
i\pa_sw -\mathcal Hw+ib\Lambda w-aw=NL(w), \ \ w=\alpha+i\beta, \ \ 
\ee
where we introduced the complex notation $w=\alpha+i\beta$, the generator of the scaling symmetry $\Lambda f=y\cdot\nabla f$ and the modulation parameters $$b=-\lsl, \ \ a=-\Theta_s,$$ and where the linearized operator is explicitely given by $$\mathcal Hw=-\Delta w+\frac{y^4-6y^2+1}{y^2(1+y^2)^2}.$$

{\bf step 2} Construction of the approximate profile and formal derivation of the law.\\
We now proceed as in \cite{MR4}, \cite{RS}, \cite{RaphRod} and look for a suitable approximate solution to the renormalized equation \fref{renormalizeequation} in the form of a homogeneous expansion $$w_0(s,y)=\alpha_0(s,y)+i\beta_0(s,y)$$ with $$\alpha_0=aT_{1,0}+b^2T_{0,2}, \ \ \beta_0=bT_{0,1}+abT_{1,1}+b^3T_{0,3}, \ \ \gamma_0=b^2S_{0,2}.$$ The goal is to {\it find the law for the modulation parameters $s\mapsto (a,b)$ allowing us to construct profiles
$T_{i,j}$} with tempered growth at infinity. At the order $b$, we get $$\mathcal HT_{0,1}=\Lambda \phi, \ \ \phi(y)=2\tan^{-1}(\frac1y)$$ which yields a growing solution for $y$ large 
$T_{0,1}(y)\sim y\log y-y\ \ \mbox{as} \ y\to +\infty.$
An explicit computation then reveals that to a leading order 
\be
\label{vnovhoe}
b_s\sim -b^2-a^2, \ \ a_s\sim 0
\ee is the unique choice which allows us to solve the $T_{i,j}$ system with controlled growth as $1\ll y$. In fact, similar to \cite{MR4}, \cite{RaphRod},  a {\it flux computation} based on the asymptotic behavior of the radiative terms $T_{i,j}$ allows us to compute the additional logarithmic corrections induced by non-trivial boundary terms at infinity: 
\be
\label{cnoehohe}
b_s+b^2\sim-\frac{b^2}{2|\log b|}-a^2, \ \ a_s\sim -2\frac{ab}{|\log b|}.
\ee A new phenomenon here is that the acceleration of the phase $\Theta$ acts as {\it a damping force against concentration} through the $b$ equation \fref{cnoehohe}, and the dynamical system \fref{cnoehohe} admits one dimensional set of initial data for which: $|a|\ll \frac{b}{|\log b|}.$ The integration of the modulation equation in this regime:
\be
\label{modulationequaitons}
b_s+b^2=-\frac{b^2}{2|\log b|}, \ \  |a|\ll\frac{b}{|\log b|}, \ \ b=-\lsl, \ \ \frac{ds}{dt}=\frac{1}{\lambda^2}, \ \ \Theta_s=-a,
\ee now yields finite time blow up $\lambda(t)\to 0\ \ \mbox{as} \ \ t\to T\ \ \mbox{for some finite}\ \ T<+\infty$ together with the asymptotics (\ref{cneoheoe}) near blow up time  and the convergence $\Theta(t)\to \Theta^*\ \ \mbox{as} \ \ t\to T$.\\

{\bf step 3} Control of the remainder: the mixed energy/Morawetz Lyapunov functional.\\
After the approximate solution is constructed, we use modulation theory to introduce a suitable nonlinear decomposition of the flow $$u(t,x)=e^{\Theta(t)R}\left[(\alpha_0+\alpha)(s,r)e_r+(\beta_0+\beta)(s,y)e_\tau+(\gamma_0+\gamma)Q\right](s,y)$$ where the four modulation parameters $(\lambda,b,\Theta,a)$ are chosen, by a standard modulation argument, to ensure that $w=\alpha+i\beta$ is orthogonal to the kernel of $\mathcal H^2$. Recall that $\mathcal H$ is a positive operator with a resonance $\mathcal H(\Lambda \phi)=0$ generated by the scaling and phase invariances. Similar to \cite{RaphRod}, our strategy to control the remainder term $w$ is to construct a Lyapunov functional mixing the energy and Morawetz type identities. The are three main differences with the analysis in \cite{RaphRod}. First we need to take more derivatives of the equation to overcome the growth of the radiation, and the Schr\"odinger map problem is, in some sense, two derivatives ``above" the wave map problem. Second, the {\it quasilinear} structure of the problem needs to be addressed through the use of suitable derivatives compatible with the geometry of the system. Third, we need to construct a codimension one set of initial data to excite the suitable solution to \fref{cnoehohe}. 
The Lyapunov type functional is built at the level of the Sobolev $H^4$ norm. Its properties in particular require the use of 
a factorization of the operator $\mathcal H$ and, thanks to the construction of a sufficiently high order approximate profile and the {\it four} orthogonality conditions on $w$, 
yield a uniform bound: 
\be
\label{vnoioueore}
\|w\|_{H^4}^2\lesssim \|\mathcal H^2 w\|_{L^2}^2\lesssim \frac{b^4}{|\log b|^2}
\ee 
Such an estimate is sufficient to control the error terms arising in the problem and, in particular, to 
verify the modulation equations \fref{modulationequaitons}.

\begin{acknowledgement}
P.R is supported by ANR Jeune Chercheur SWAP. I. R. is supported by NSF grant DMS-1001500.
\end{acknowledgement}


\begin{thebibliography}{10}
\bibitem{BT} Bejenaru, I.; Tataru, D.;
\newblock Near soliton evolution for equivariant Schr\"odinger Maps in two spatial dimensions;
\newblock arXiv:1009.1608.
%
\bibitem{BT2} Bejenaru, I.; Ionescu, A.;  Kenig, C.;  Tataru, D.;
 \newblock Global Schr\"odinger maps;
 \newblock to appear, Annals of Math..

\bibitem{heatflow} Van den Bergh, J.; Hulshof, J.; King, J., Formal asymptotics of bubbling in the harmonic map heat flow, SIAM J. Appl. Math. vol 63, o5. pp 1682-1717.

\bibitem{CSU} Chang, N-H.; Shatah, J.; Uhlenbeck, K. 
\newblock Schr\"odinger maps. 
\newblock Comm. Pure Appl. Math. 53 (2000), no. 5.

\bibitem{GS} Grotowski, J.; Shatah, J. 
\newblock Geometric evolution equations in critical dimensions. 
\newblock Calc. Var. Partial Differential Equations 30 (2007), no. 4, 499-512.

\bibitem{T3} Gustafson, S.; Kang, K.; Tsai, T-P.;
\newblock Schr\"odinger flow near harmonic maps;
\newblock Comm. Pure Appl. Math. 60 (2007), no. 4, 463-499. 

\bibitem{T4} Gustafson, S.; Kang, K.; Tsai, T-P.;
\newblock Asymptotic stability of harmonic maps under the Schr\"odinger flow.;
\newblock Duke Math. J. 145 no. 3 (2008) 537-583.

\bibitem{T2} Gustafson, S.; Nakanishi, K.; Tsai, T-P.;
\newblock  Asymptotic stability, concentration and oscillations in harmonic map heat flow, Landau Lifschitz and Schr\"odinger maps on $\Bbb R^2$;
\newblock Comm. Math. Phys. (2010), 300, no 1, 205-242.

%
%


%
\bibitem{KS} Krieger, J.; Schlag, W.;
\newblock Concentration compactness for critical wave maps;
\newblock Preprint, 2009.

\bibitem{KST} Krieger, J.; Schlag, W.; Tataru, D. ;
\newblock Renormalization and blow up for charge one equivariant critical wave maps;
\newblock Invent. Math. 171 (2008), no. 3, 543-615. 


\bibitem{MR1} Merle, F.; Rapha\"el, P.; 
\newblock Blow up dynamics and upper bound on the blow up rate for critical nonlinear Schr\"odinger equation, 
\newblock Annals of Math. {\bf 161} (2005), no. 1, 157-222.
%
\bibitem{MR2} Merle, F.; Rapha\"el, P.;
\newblock  Sharp upper bound on the blow up rate for critical nonlinear Schr\"odinger equation,
\newblock Geom. Funct. Anal. 13 (2003), 591-642.

%
\bibitem{MR3} Merle, F.; Rapha\"el, P.;
\newblock On universality of blow up profile for $L^2$ critical nonlinear Schr\"odinger equation,
\newblock  Invent. Math. 156, 565-672 (2004).

%
\bibitem{MR4} Merle, F.; Rapha\"el, P.;
\newblock Sharp lower bound on the blow up rate for critical nonlinear Schr\"odinger equation,
\newblock  J. Amer. Math. Soc. 19 (2006), no. 1, 37-90.

%

%
\bibitem{MR5} Merle, F.; Rapha\"el, P.;
 \newblock Profiles and quantization of the blow up mass for critical non linear Schr\"odinger equation, 
 \newblock Comm. Math. Phys. {\bf 253} (2004), no. 3, 675-704. 


\bibitem{MRR} 
Merle, F.; Rapha\"el, P.; Rodnianski, I.;
\newblock Blow up for the energy critical corotational Schr\"odinger map
\newblock Preprint 2011.

%
\bibitem{Q} Qing, J.;
 \newblock On singularities of the heat flow for harmonic maps from surfaces into spheres, 
 \newblock Comm. Anal. Geom. 3 (1995), 297-315.
 
\bibitem{QT} Qing J.;  Tian G.;
 \newblock Bubbling of the heat flows for harmonic maps from surface, 
 \newblock Comm. Pure Appl. Math. 50 (1997),  295-310.


%
\bibitem{R1} Rapha\"el, P.;
\newblock Stability of the log-log bound for blow up solutions to the critical nonlinear Schr\"odinger equation
\newblock Math. Ann. 331 (2005), 577-609.
%
\bibitem{RaphRod} Rapha\"el, P.; Rodnianski. I.;
 \newblock Stable blow up dynamics for the critical co-rotational Wave map and equivariant Yang Mills problems,
 \newblock Preprint 2009.
 
 \bibitem{RS} Rapha\"el, P.; Szeftel, J.;
 \newblock Existence and uniqueness of minimal blow up solutions to an inhomogeneous mass critical NLS equation;
 \newblock to appear in Jour. Amer. Math. Soc..
 
 
 \bibitem{RSt} 
 \newblock Rodnianski, I.; Sterbenz, J. 
 \newblock On the formation of singularities in the critical $O(3)$ $\sigma$-model;
 \newblock Ann. of Math. (2) 172 (2010), no. 1, 187-242.
 
 \bibitem{PS}
 \newblock Smith, P.;
 \newblock Conditional global regularity of Schr\"odinger maps: sub-threshold dispersed energy;
 \newblock Preprint 2010.
 
 \bibitem{ST} Sterbenz, J.; Tataru, D.;
 \newblock Regularity of Wave-Maps in dimension 2 + 1;
 \newblock Comm. Math. Phys. {\bf 298} (2009), 139-230.
 
\bibitem{S} Struwe, M.;
\newblock
On the evolution of harmonic mappings of Riemannian surfaces, 
\newblock Comm. Math. Helv. 60, 558-581 (1985).

\bibitem{T} Tao, T.; 
\newblock Global regularity of wave maps III-VII,
\newblock Preprints 2008-2009.

\bibitem{PT} Topping, P.M.;
\newblock Winding behaviour of finite-time singularities of the harmonic map heat flow;
\newblock Math. Z. 247 (2004).

\end{thebibliography}
\end{document}